\documentclass{amsart}

\usepackage{amsfonts,amssymb,amscd,amsmath,enumerate,verbatim,calc,xypic}
\usepackage[all]{xy} \SelectTips{eu}{}
\def\wh{\widehat}

\newcommand{\bsx}{{\boldsymbol x}}

\def\CX{{\mathcal X}}

\newcommand{\BN}{{\mathbb N}}

\newcommand{\fm}{{\mathfrak m}}
\newcommand{\fn}{{\mathfrak n}}

\newcommand{\bsf}{{\boldsymbol f}}

\newcommand{\im}{\operatorname{Im}}

\newcommand{\Soc}{\operatorname {Soc}}

\newcommand{\col}{\colon}

\newcommand{\Ker}{\operatorname{Ker}}

\newcommand{\pd}[2]{\operatorname{pd}_{#1}{#2}}
\newcommand{\ci}[2]{\operatorname{CI-dim}_{#1}{#2}}

\newcommand{\Hom}[3]{\operatorname{Hom}_{#1}({#2},{#3})}

\newcommand{\Ext}[4]{\operatorname{Ext}^{#1}_{#2}(#3,#4)}
\newcommand{\Tor}[4]{\operatorname{Tor}_{#1}^{#2}(#3,#4)}

\newcommand{\edim}{\operatorname{edim}}

\newcommand{\codim}{\operatorname{codim}}

\newcommand{\curv}[2]{\operatorname{curv}_{#1}{(#2)}}

\newcommand{\cx}[2]{\operatorname{cx}_{#1}({#2)}}
\newcommand{\px}[2]{\operatorname{px}_{#1}({#2)}}
\newcommand{\cxx}[3]{\operatorname{cx}_{#1}({#2},{#3})}

\theoremstyle{plain}
\newtheorem{theorem}{Theorem}[section]
\newtheorem{corollary}[theorem]{Corollary}
\newtheorem{lemma}[theorem]{Lemma}

\newtheorem{question}[theorem]{Question}

\newtheorem{proposition}[theorem]{Proposition}

\theoremstyle{definition}

\newtheorem{definition}[theorem]{Definition}
\newtheorem{chunk}[theorem]{}

\newtheorem{remark}[theorem]{Remark}

\numberwithin{equation}{theorem}

\theoremstyle{remark}

\newcommand{\numberseries}{\mdseries}   %Fontseries used for numbering theorem

\newlength{\thmtopspace}                %Space above theorem
\newlength{\thmbotspace}                %Space below theorem
\newlength{\thmheadspace}               %Space between theorem caption and text
\newlength{\thmindent}                  %For indenting

\newtheoremstyle{bfupright head,italic body}
                {\thmtopspace}{\thmbotspace}
                {\itshape}{\thmindent}{\bfseries}{}{\thmheadspace}
                {{\numberseries \thmnumber{\bf(#2)}}\thmnote{#3}}

\theoremstyle{bfupright head,italic body}
             \newtheorem*{ires*}{}

\begin{document}

\date{\today}

\title{Comparing complexities of pairs of modules}

\dedicatory{Dedicated to Professor Paul C. Roberts on the occasion of
his sixtieth birthday}

\author[H.~Dao]{Hailong~Dao}

\address{Department of Mathematics, University of Kansas,  405 Snow Hall, 1460 Jayhawk Blvd, Lawrence,
KS 66045-7523, USA} \email{hdao@math.ku.edu}

\author[O.~Veliche]{Oana~Veliche}
\address{Department of Mathematics, University of Utah, Salt Lake City, Utah~84112}
\email{oveliche@math.utah.edu}

%\thanks
\subjclass [2000]{13D07, 13H10}

\keywords{Complexity, complexity inequalities, artinian rings,
Cohen-Macaulay isolated singularity, Auslander-Reiten}
\begin{abstract}
Let $R$ be a local ring and $M,N$ be finitely generated $R$-modules.
The complexity of $(M,N)$, denoted by $\cxx RMN$, measures the
polynomial growth rate of the number of generators of the modules
$\Ext nRMN$. In this paper we study several basic equalities and
inequalities involving complexities of different pairs of modules.
\end{abstract}
\maketitle
%\tableofcontents

\section{Introduction}
\label{Introduction}

Let  $R$ be  a commutative local noetherian ring with maximal
ideal $\fm$  and residue field $k=R/\fm$, and let $M$ and $N$ be
finitely generated $R$-modules.
The {\it complexity} of the pair of modules $(M,N)$,
denoted by $\cxx{R}MN$, measures the polynomial growth rate of the
number of generators of the modules $\Ext nRMN$; see Section 2 for
background and definitions. It was first introduced  by Avramov and
Buchweitz in \cite{ab00} to study properties of $\Ext nRMN$  when
$R$ is a complete intersection. Over such rings, properties of the
complexity of a pair of modules have been studied extensively;
see e.g. \cite{ab00,agp97,d}. In this paper we study the complexity
of a pair of modules over rings other than complete intersections.

Avramov and Buchweitz prove in \cite{ab00} that when
$R$ is a complete intersection,  $\cxx RMN$ cannot exceed
either $\cx RM\col\!\!\!=\cxx RMk$, the complexity of $M$, or $\px RN\col\!\!\!= \cxx
RkN$, the plexity of $N$. Thus, we ask the following

\begin{question}
\label{q1} Let $R$ be a local noetherian ring. Is it true that the
inequality
\begin{equation*}
\cxx RMN\leq\min\{\cx RM,\px RN\}
\end{equation*}
holds for all finitely generated $R$-modules $M$ and $N$?
\end{question}
\noindent  Note that if the right-hand side is zero then the
left-hand side is automatically zero.
We show that an affirmative answer holds for artinian rings,
see Lemma \ref{cx-px-inequalities},
and more generally for local Cohen-Macaulay rings with isolated
singularity, see Theorem \ref{cx-inequality}.

Another motivation for our study is a number of questions related to
the Auslander-Reiten Conjecture, which asserts that over a local
ring $R$, a module $M$ with $\Ext iRM{M\oplus R}=0$ for all $i>0$
must be free. To highlight the connection with complexity, we first
formulate an asymptotic version of this conjecture, which has
implicitly appeared in some recent papers, see Remark
\ref{AAR-rings}. We say that a ring $R$ has the {\it asymptotic
Auslander-Reiten property} if it satisfies:

\noindent (AAR)
\hspace{0.5cm} For any finitely generated $R$-module $M$ the implication
\begin{equation*}
\cxx RMR=0=\cxx RMM \implies \cx RM=0
\end{equation*}
\hspace{1.7cm} holds.

A ring with (AAR) property satisfies the Auslander-Reiten conjecture
by Remark \ref{AR and AAR}. In this paper, we focus on the following
properties of a ring $R$, which are stronger then (AAR):
\begin{equation*} \tag{P1} \label{p1} \cxx RMR=\cx
RM\quad\text{for all finitely generated $R$-modules}\  M.
\end{equation*}
\begin{equation*}
\tag{P2} \label{p2} \cxx RMM=\cx RM\quad\text{for all finitely
generated $R$-modules}\  M.
\end{equation*}

Our investigation identifies certain classes of local, artinian
rings satisfying the properties described above. For example, an
artinian ring $(R,\fm)$ satisfies the property \eqref{p1} if
$2\ell_R(\Soc(R))> \ell_R(R)$ or if $R$ is non-Gorenstein with
$\fm^3=0$ and $2\ell_R(\Soc(R))> \ell_R(R)-2$; see Proposition
\ref{2r>l} and Theorem \ref{cx m3}. On the other hand, Gorenstein
rings with radical cube zero satisfy  \eqref{p2}; see Proposition
\ref{cx m3-Gorenstein}. Note that complete intersection rings
satisfy \eqref{p2}; see \cite[Theorem II]{ab00}. One interesting
feature of our results is that non-regular rings satisfying
\eqref{p1} are far from being Gorenstein, while the ones satisfying
\eqref{p2} form a strict subclass of artinian Gorenstein rings; see
Remarks \ref{remark Gorenstein} and \ref{LJ}.

The structure of the paper is summarized below. Section
\ref{Preliminaries} describes some preliminary results. In Section
\ref{Complexity of modules over artinian rings} we prove several
inequalities and equalities of complexities over artinian rings.
Section \ref{Rings with isolated singularity} is devoted to the
proof of Theorem \ref{cx-inequality} which asserts that over a
Cohen-Macaulay local ring $R$ with isolated singularity, the
inequality $\cxx RMN\leq\min\{\cx RM,\px RN\}$ holds. In Sections
\ref{AAR} we study rings satisfying properties \eqref{p1} and
\eqref{p2}.

%%%%%%%%%%%%%%%%%%%%%%%%% SECTION 1 %%%%%%%%%%%%%%%%%%%%%%%%%%%%%%%%%%%%%%%%%%

\section{Preliminaries}
\label{Preliminaries}

In this section,  we recall the definition of the complexity of a
sequence  and of a pair of modules, and then prove and recall some
of their properties used throughout the paper.

\begin{definition}
The  {\it complexity of the sequence} $\{x_i\}_{i\geq0}$, of non-negative numbers is given by
$$\cx {}{\{x_i\}} = \inf\left\{b\in\BN\Bigg|\begin{gathered}x_i\leq a\cdot i^{b-1}\ \text{for some}\\ \text{real number}\ a\  \text{and for all}\ i\gg 0\end{gathered}\right\}.$$
\end{definition}

\begin{proposition}
\label{sequences} Let $\{x_i\}_{i\geq0}$ and $\{y_i\}_{i\geq0}$ be
sequences of non-negative integers. Let $a,b$ be positive real
numbers.
\begin{enumerate}[\rm\quad(1)]
\item If  $ a\cdot y_i \leq x_i \leq b\cdot y_i$ for  all $i\gg 0,$
then $\cx{}{\{y_i\}} = \cx{}{\{x_i\}}.$

\item $\cx{}{\{x_{i+1}-x_i\}}\geq \cx{}{\{x_i\}}-1.$

\item If $y_i=a\cdot x_{i+1}+b\cdot x_{i}$, then $\cx{}{\{y_i\}}=\cx{}{\{x_i\}}$.
\item If $y_i = a\cdot x_{i+1}-b\cdot x_{i}$ and $a>b$, then $\cx{}{\{y_i\}}=\cx{}{\{x_i\}}.$
\end{enumerate}
\end{proposition}

\begin{proof}
The proofs of (1), (2), (3) are straightforward.

(4): Set $d = \cx{}{\{x_i\}}$. Since $\cx{}{\{y_i\}} \leq
\cx{}{\{ax_{i+1}\}}=d$ it is enough  to prove $d\leq
\cx{}{\{y_i\}}$. We consider the following three cases.

Case $d=0$ is trivial.

Assume $d$=1. Suppose $\cx{}{\{y_i\}}=0$. Then $a\cdot
x_{i+1}-b\cdot x_{i}=0$ for $i\gg0$. But since  $\{x_i\}$ is a
bounded sequence of integers and $a>b$, we must have $x_i=0$ for
$i\gg 0$, a contradiction.

Assume  $d\geq 2$. There exists a positive integer $L$ such that
$x_{L} \geq x_{L-1}$ and, by definition, there exists a subsequence
$\{x_{i_j}\}_{j\geq 0}$ such that $i_0\geq L$ and
\begin{equation*}
\tag{$*$} \lim_{j \to \infty} {x_{i_j}}/{(i_j)^{d-2}} =\infty.
\end{equation*}

For each $j\geq 0$, let $t_j$ be the biggest integer such that
\begin{equation*}
\tag{$**$} t_j \leq i_j  \quad \text{and} \quad x_{t_j} \geq
x_{t_{j}-1}.
\end{equation*}

Note that such $t_j$ exists because $i_j\geq L$ for all $j\geq 0$.
Our choice of $t_j$ ensures that $x_{t_j} \geq x_{i_j}\quad\text{for
all}\ j\geq 0$. This together with $(*)$ gives
$$\lim_{j \to \infty} {x_{t_j}}/{(t_j)^{d-2}} =\infty.$$
On the other hand, by $(**)$ we get $ y_{t_j-1} \geq (a-b)x_{t_j}
\geq x_{t_j}.$ Thus, $$\lim_{j \to \infty} y_{t_j-1}/{(t_j-1)}^{d-2}
=\infty.$$ This implies that $\cx{}{\{y_i\}} \geq d$, which is what
we need.
\end{proof}

For the rest of this section, let $R$ be a commutative local
noetherian ring with maximal ideal $\fm$  and residue field
$k=R/\fm$, and  let $M,N$ be finitely generated $R$-modules.

\begin{definition}
\label{defcx}
In \cite{ab00} Avramov and Buchweitz define the {\it complexity of the pair of modules} $(M,N)$ to be
$$\cxx RMN=\cx{}{\{\nu_R(\Ext iRMN)\}},$$
where $\nu_R(-)$ denotes the minimal number of generators. Set
$$\cx RM=\cxx RMk\quad\text{and}\quad\px RM=\cxx RkM.$$
\end{definition}

\noindent
It is easy to see that
$$\cxx RMN=0\quad\text{if and only if}\quad\Ext iRMN=0\quad\text{for all}\ i\gg 0.$$
Another immediate property is the following.

\chunk
\label{cx-regular}
If $x \in R$ is an $R$-regular
and $M$-regular element such that $xN=0$, then
$$\cxx RMN=\cxx {R/xR}{M/xM}N.$$

In general, it is easier to work with the length function
$\ell_R(-)$, if possible, than with the function $\nu_R(-)$; the
former is additive on short exact sequence while the latter is not.
Thus, the following easy result will be very useful.

\begin{lemma}
\label{mh-vanishing}
Let $R$ be a local ring and let $\{M_i\}_{i\geq
0}$ be  a sequence of finitely generated $R$-modules. Suppose there
exists a positive integer $h$ such that $\fm^hM_i=0\ \text{for all}\
i\gg 0$. Then
$$\cx {}{\{\nu_R(M_i)\}}=\cx{}{\{\ell_R(M_i)\}}.$$
\end{lemma}

\begin{proof}
We may assume that $\fm^hM_i=0\ \text{for all}\  i\geq 0$, thus  $M_i$ is an $R/\fm^h$ module for all $i\geq 0$. Therefore, we
get the inequalities
$$ \frac{1}{\ell_R(R/\fm^h)}\ell_R{(M_i)}\leq\nu_R{(M_i)}\leq\ell_R{(M_i)}\quad\text{for all}\ i\geq 0.$$
The conclusion follows by Proposition \ref{sequences}(1).
\end{proof}

\begin{corollary}
\label{nu-ell}
If $R$ is an artinian local ring, then in the
definition of $\cxx RMN$, one can replace the function $\nu_R(-)$ by
the length function $\ell_R(-)$.
\end{corollary}

\begin{proof}
We may assume that $M$ is non-zero. Since $R$ is artinian, there exists a
positive integer $h$ such that $\fm^hM=0$. In particular, we have
$\fm^h\Ext iRMN=0$ for all $i\geq 0$; now  we apply  Lemma
\ref{mh-vanishing}.
\end{proof}

\begin{remark}
Over  an artinian ring, whenever we work with the complexity of a
pair of  finitely generated $R$-modules, we can use Corollary
\ref{nu-ell} and work with the length function.
\end{remark}

Finally, we recall some known results on complexity that we refer to
in the paper.

\chunk
\label{k ci}
\cite[(8.1.2)]{a98}
The local ring $(R,\fm, k)$
is a complete intersection if and only if $\cx Rk <\infty$.

\chunk
\cite[Theorem II]{ab00}
\label{properties-ci}
If $R$ is a
local complete intersection and $M,N$ are finitely generated
$R$-modules, then
\begin{enumerate}[\rm\quad(1)]
\item $\cxx RMM=\cx RM=\px RM <\infty$.
\item $\cx RM+\cx RN-\codim R\leq\cxx RMN=\cxx RNM.$
\item $\cxx RMN\leq\min\{\cx RM,\cx RN\}$.
\end{enumerate}

\chunk \label{E Gorenstein}\cite[(1.1)]{jl07} Let $(R,\fm,k)$ be a
local ring such that $\fm^3=0$ and let $E$ be the injective
envelope of $k$. Then $R$ is Gorenstein if and only if $\cx
RE<\infty$.

%%%%%%%%%%%%%%%%%%%%%%%%%%%%%%%%%%%%%%%%%%% SECTION 2 %%%%%%%%%%%%%%%%%%%%%%%%%%%%%%%

\section{Complexity of modules over artinian rings}
\label{Complexity of modules over artinian rings}

In this section, $R$ is a local artinian ring and  $M,N$ are
finitely generated $R$-modules. We study basic inequalities and
equalities related to Question \ref{q1} from the introduction. A
technical but useful result is Lemma \ref{N-inequality} which
establishes an inequality between the length of the modules $\Ext
nRMN$ and certain Betti numbers  of $M$ and $N$.

\chunk
\label{cx-px}
Let  $E$  be  the injective envelope of the
residue field $k$ and let $M^\vee=\Hom RME$ be the Matlis dual of
$M$. There are isomorphisms ${\Ext iRMN}^\vee\cong\Tor iRM{N^\vee}$
for all $i>0$. In particular, we get the equalities
\begin{equation*}
\cxx RMN=\cxx R{N^\vee}{M^\vee},\quad\text{and}\quad\cx RM=\px
R{M^\vee}.
\end{equation*}

\begin{lemma}
\label{cx-px-inequalities}
If  $R$ is an artinian local ring, then for every
finitely generated $R$-modules $M$ and $N$, we have the inequality
\begin{equation*}
\cxx RMN\le\min\{\cx RM,\px RN\}.
\end{equation*}
\end{lemma}

\begin{proof}
Set $b_i=\beta^R_i(M)$, the $i$-th Betti number of $M$  for all
$i\geq 0$ and consider a minimal free resolution of the module $M$
$$\cdots\to R^{b_{i+1}}\xrightarrow{\partial_{i+1}}R^{b_i}\xrightarrow{\partial_i} R^{b_{i-1}}\to\cdots\to R^{b_1}\xrightarrow{\partial_1} R^{b_0}\xrightarrow{\partial_0}0\to\cdots.$$
Applying the functor $\Hom R{-}N$, we get the complex
$$\cdots\to N^{b_{i-1}}\xrightarrow{\Hom R{\partial_{i}}N}N^{b_i}\xrightarrow{\Hom R{\partial_{i+1}}N} N^{b_{i+1}}\to\cdots.$$
By definition $\Ext iRMN$ is a homomorphic image of $\Ker(\Hom
R{\partial_{i+1}}N)$, so
$$\ell_R(\Ext iRMN)\leq \ell_R(\Ker(\Hom R{\partial_{i+1}}N))\leq \ell_R(N^{b_i})=b_i\ell_R(N).$$
Thus, the  inequality $\cxx RMN\leq \cx RM$ holds because of
\ref{nu-ell}; the inequality $\cxx RMN\leq\px RN$ follows from the
first one  and \ref{cx-px}.
\end{proof}

\begin{lemma}
\label{N-inequality} Let $(R,\fm)$ be an artinian local ring and let
$M$ and $N$ be finitely generated $R$-modules. If $I$ is an
ideal of $R$ such that $(I\fm)N=0$, then for every $i\geq 0$ we have
the inequality
$$\ell_R(\Ext iRMN) \geq
\ell_R(N)\cdot\beta_i^R(M)-\ell_R(R/I)\nu_R(N)\cdot[\beta_{i-1}^R(M)
+ \beta_i^R(M)].$$

\end{lemma}

\begin{proof}
As above, set $b_i=\beta^R_i(M)$ for all $i\geq 0$ and consider a minimal free resolution of the module $M$ with differential $\partial=\{\partial_i\}_{i\geq 0}$. Set
$$K_i = \Ker(\Hom R{\partial_{i+1}}N)\quad\text{and}\quad C_i=\im(\Hom R{\partial_{i}}N).$$
By definition we have  $ \Ext iRMN = K_i/C_i$ and the exact sequence
\begin{equation}
\tag{$*$} 0\to K_i\to N^{b_i}\to C_{i+1}\to 0.
\end{equation}
Thus, we obtain the equalities
\begin{align*}
\ell_R(\Ext iRMN)&= \ell_R(K_{i}) - \ell_R(C_i),\quad\text{and}\\
\ell_R(N^{b_i})&= \ell_R(K_{i})+\ell_R(C_{i+1}).
\end{align*}
By elimination, we get
\begin{equation}
\tag{$**$}
\ell_R(\Ext iRMN)= \ell_R(N)\cdot b_i- \ell_R(C_{i+1}) - \ell_R(C_i),\quad\text{for all}\quad i\geq 0.
\end{equation}

Since $(I\fm)N=0$ and $C_i\subseteq \fm N^{b_i}$, we get $IC_i=0$ for all  $i\geq 0$. This implies that $ \ell_R(R/I) \nu_R(C_i)\geq \ell_R(C_i)$ for all $i\geq 0$, which together with $(**)$ implies the inequality
$$\ell_R(\Ext iRMN) \geq \ell_R(N)\cdot b_i-\ell_R(R/I)\cdot [\nu_R(C_{i+1}) + \nu_R(C_{i})]. $$
Observe that $\nu_R(C_i) \leq \nu_R(N)\cdot b_{i-1}$ by the exact
sequence ($*$). Thus, we get:
\begin{equation*}
\ell_R(\Ext iRMN) \geq\ell_R(N)\cdot b_i-\ell_R(R/I)\nu_R(N)\cdot (b_{i} + b_{i-1}).
\end{equation*}
which is what we want.
\end{proof}

\begin{proposition}
\label{cx MN} Let $(R,\fm)$ be a local artinian ring and let $M$ and
$N$ be finitely generated $R$-modules such that $(I\fm)N=0$ for
some ideal $I$  of $R$.
\begin{enumerate}[\rm\quad(1)]
\item If $\ell_R(N)>2\ell_R(R/I)\nu_R(N)$, then
$\cxx RMN = \cx RM.$\\
If  $\ell_R(N)=2\ell_R(R/I)\nu_R(N)$, then $\cxx RMN\in\{\cx RM-1,\ \cx RM\}.$
\item If $\ell_R(N) > 2\ell_R(R/I)\nu_R(N^\vee),$  then
$\cxx RNM = \px RM.$\\
If  $\ell_R(N)= 2\ell_R(R/I)\nu_R(N^\vee)$, then $\cxx RNM\in\{\px RM-1,\ \px RM\}.$
\end{enumerate}
 \end{proposition}

\begin{proof}
(1): Set $a=\ell_R(N)-\ell_R(R/I)\nu_R(N)$ and
$b=\ell_R(R/I)\nu_R(N)$. Then  by Lemma \ref{N-inequality} we have
the inequality $$\ell_R(\Ext iRMN) \geq a\cdot \beta_i^R(M) -b\cdot
\beta_{i-1}^R(M).
$$ The conclusions follow by  Lemma \ref{cx-px-inequalities} and Proposition
\ref{sequences}.\\
(2) We apply part (1) to the pair of modules $(M^\vee,N^\vee)$. One
needs to use the facts that $\ell(N) = \ell(N^{\vee})$,
$(I\fm)N^\vee=0$, and the equalities in \ref{cx-px}.
\end{proof}

\begin{proposition}
\label{m2-cx} Let $(R,\fm)$  be a local artinian  ring and let $M$
and $N$ be finitely generated $R$-modules with $\fm^2N=0$.
\begin{enumerate}[\rm\quad(1)]
\item If $\ell_R(\fm N)>\nu_R(N)$, then $\cxx RMN=\cx RM$.
\item If $\ell_R(\fm N)<\nu_R(N)$, then $\cxx RNM=\px RM$.
\item If $\ell_R(\fm N)=\nu_R(N)$, then
\begin{align*}
\cxx RMN&\in\{\cx RM-1,\ \cx RM\}\quad\text{and} \\
\cxx RNM &\in\{\px RM-1,\ \px RM\}.
\end{align*}
\end{enumerate}
\end{proposition}

\begin{proof}
(1) and the first inclusion of (3) follow directly from Proposition \ref{cx MN}(1).

(2) and the second inclusion of (3): Set $a=\nu_R(N)$ and
$b=\ell_R(\fm N)$. There is  an exact sequence $0\to k^b\to N\to
k^a\to 0$. Applying  $\Hom R-M$ we get a long exact sequence
$$\cdots\to k^{b\cdot\mu_R^i(M)}\to k^{a\cdot\mu_R^{i+1}(M)}\to \Ext {i+1}RNM\to k^{b\cdot\mu_R^{i+1}(M)}\to\cdots$$
Using the additivity of length  we get the inequalities
\begin{align*}
\ell_R(\Ext{i+1}RNM)&\geq a\cdot\mu_R^{i+1}(M)-b\cdot\mu_R^i(M)
\end{align*}
The conclusions now follow from Proposition \ref{sequences}.
\end{proof}

\begin{corollary}
\label{m2N} Let $(R,\fm)$  be a local artinian  ring  and let $N$ be
a finitely generated $R$-module with $\fm^2N=0$. Then either $R$ is
a complete intersection, and
$$\cx RN=\px RN=\cxx RNN<\infty,$$
or $R$ is not a complete intersection and one of the following
(possibly both) holds:
\begin{enumerate}[\rm\quad(1)]
\item $\px RN=\infty\quad\text{and}\quad\cx RN=\cxx RNN.$
\item $\cx RN=\infty\quad\text{and}\quad\px RN=\cxx RNN.$
\end{enumerate}
In particular, if $\cxx RNN=0$, then the module $N$ is free or injective.
\end{corollary}

\begin{proof}
If $R$ is a complete intersection, then one can apply
\ref{properties-ci}(1).

Assume that $R$ is not a complete intersection, thus $\cx Rk=\px
Rk=\infty$; see  \ref{k ci}. If $\ell_R(\fm N)>\nu_R(N)$, then
applying Proposition \ref{m2-cx}(1) to the pairs of modules $(k,N)$
and $(N,N)$ gives us
 case (1). If $\nu_R(N)>\ell_R(\fm N)$, then we apply Proposition \ref{m2-cx}(2)  to  $(k,N)$ and $(N,N)$ to get case (2). Finally, if
$\ell_R(\fm N)=\nu_R(N)$, then we apply Proposition
\ref{m2-cx}(3) to the pairs $(k,N)$ and $(N,N)$. In this situation
we get $\px RN=\cx RN=\cxx RNN=\infty$, that is (1) and (2).
\end{proof}

%%%%%%%%%%%%%%%%%%%%%%%%%%%%%%%%%%%%%% SECTION 3%%%%%%%%%%%%%%%%%%%%%%%%%%%%%%%%%%%%%%%%%%%%%%%%

\section{Rings with isolated singularity}
\label{Rings with isolated singularity}

The main result of this section, whose proof is given at the end, is
the following.

\begin{theorem}
\label{cx-inequality} Let $(R,\fm)$ be a  Cohen-Macaulay local ring.
If $M$ is  a finitely generated $R$-module such that
$\pd{R_p}{M_p}<\infty$ for any prime ideal $p\not=\fm$, then for
every finitely generated $R$-module $N$ we have the inequality
$$\cxx RMN\leq \min \{\cx RM,\px RN\}.$$
\end{theorem}

The following consequence of this result gives a partial answer to
Question \ref{q1}.

\begin{corollary}
\label{cor-cx-inequality} If $R$ is a Cohen-Macaulay local ring
with isolated singularity, then for all finitely generated
$R$-modules $M$ and $N$, we have the inequality
$$\cxx RMN\leq \min \{\cx RM,\px RN\}.$$
\end{corollary}

Next, we prove a series of preparatory results.

\begin{lemma}
\label{same h} Let $(R,\fm,k)$ be a local ring. Let $M$ be a
finitely generated $R$-module such that $M_p$ is a free $R_p$-module
for any prime ideal $p\not=\fm$. There exists a positive integer $h$
such that $\fm^h\Ext iRMN=0$ for all $i>0$ and for all $R$-modules
$N$.
\end{lemma}

\begin{proof}
Set $\CX=\{x\in \fm\mid M_x\ \text{is a free}\ R_x-\text{module}\}$
and let $I$ be the ideal of $R$ generated by all elements of $\CX$.
Without loss of generality, we may assume that $M$ is not free,
so $I$ is a proper ideal.

First, we show that $I$ is an $\fm$-primary ideal. If it is not,
then there exists a prime ideal $p \neq \fm$ such that $I\subseteq
p$. As $M_p$ is free, there exists $y\in\fm\setminus p$ such that
$M_y$ is a free $R_y$-module; that is a contradiction.

Second, we claim that for each  $x\in\CX$ there is a non-negative
integer $n(x)$ such that $x^{n(x)}\Ext iRMN=0 $ for all $i>0$ and
for all $R$-modules $N$. We have $M_x \cong
R_x^s$ for some $s$. This isomorphism is induced by a homomorphism
of $R$-modules $f\colon R^s \to M$. Let $Z,B$ and $C$ be the kernel,
image and  the cokernel of $f$ respectively. Since $f$ becomes an
isomorphism after localizing at $x$, there is an integer $n(x)$ such
that $x^{n(x)}Z=0$ and $x^{n(x)}C=0$. The long exact sequences of
Ext, obtained after applying the functor $\Hom R-N$  to the short
exact sequences
\begin{align*}
0 &\to Z \to R^s \to B\to 0\quad\text{and}\\
0 &\to B \to M \to C \to 0
\end{align*}
show that $x^{n(x)}\Ext iRMN =0$ for all $i>0$.

Finally, since $R$ is noetherian we may choose a subset
$\{x_1,\dots,x_l\}$ of $\CX$ whose elements generate the ideal $I$.
Let $n=\max\{n(x_1),\dots,n(x_l)\}$. Then $I^{nl}\Ext iRMN=0$ for
all $i>0$ and for all $R$-modules $N$.
Since $I$ is $\fm$-primary the desired conclusion now
follows.
\end{proof}

\begin{lemma}
\label{x-reduction} Let $(R,\fm)$ be a local Cohen-Macaulay ring of
Krull dimension $d$. Let $M$ and $N$ be maximal Cohen-Macaulay
$R$-modules such that $M_p$ is free for all prime ideals
$p\not=\fm$. Then there exists an $R$-regular sequence $\bsx$ of
length $d$ such that
\begin{equation*}
\cxx RMN=\cxx {R/\bsx R}{M/\bsx M}{N/\bsx N}.
\end{equation*}
\end{lemma}

\begin{proof}
We construct the sequence $\bsx$  of length $d$ inductively.

If $d=0$ there is nothing to be proved.

Assume $d\geq 1$. In this case it is enough to find the first
element $x_1$ of the sequence. Indeed, the modules $M/x_1M$ and
$N/x_1N$ are maximal Cohen-Macaulay and  $M/x_1M$ is free  on the
punctured spectrum of $R/x_1R$, thus we can continue inductively.

By  Lemma \ref{same h} there exists a positive integer $h$ such that
$\fm^h\Ext iRMN=0$ for all $i>0$  and all
$R$-modules $N$. Choose an
$R$-regular element $x_1$ in $\fm^h$. Since $M$ and $N$ are maximal
Cohen-Macaulay, $x_1$ is also $M$- and $N$-regular.

First, we show that
\begin{equation*}
\tag{$*$}
\cxx RMN=\cxx RM{N/x_1N}.
\end{equation*}
From the short exact sequence $0\to N\xrightarrow{ x_1} N\to
N/x_1N\to 0$ we obtain by applying the functor $\Hom RM-$, the long
exact sequence
$$\cdots\to\Ext iRMN\xrightarrow{ x_1}\Ext iRMN\to\Ext iRM{N/x_1N}\to\Ext{i+1}RMN\to\cdots.$$
Since $x_1$ is in $\fm^h$,  this long exact sequence splits into short exact sequences of modules of finite length
$$0\to\Ext iRMN\to\Ext iRM{N/x_1N}\to\Ext {i+1}RMN\to 0\quad\text{for all}\ i>0,$$
as  $\fm^h\Ext iRM{N/x_1N}=0$ for all $i>0$.

Using the additivity of the length function we get for all $i>0$
$$\ell_R(\Ext iRM{N/x_1N})=\ell_R(\Ext iRMN)+\ell_R(\Ext {i+1}RMN).$$
Applying now Lemma \ref{mh-vanishing} and Proposition \ref{sequences}(3) we obtain the desired equality of complexities.

Second, by \ref{cx-regular} we have
\begin{equation*}
\tag{$**$}
\cxx RM{N/x_1N}=\cxx{R/x_1R}{M/x_1M}{N/x_1N}.
\end{equation*}

Combining now the equalities $(*)$ and $(**)$ finishes the proof.
\end{proof}

\begin{lemma}
\label{plex_deform}
Let $R$ be a local ring, $N$ a finitely generated  $R$-module, and let $x$ be a regular
element on $R$ and $N$. Then, for any finitely generated  $R/xR$-module  $M$ we have the equality
$$\cxx RMN=\cxx {R/xR}M{N/xN}.$$
In particular, $\px RN = \px {R/xR}{N/xN}$.
\end{lemma}

\begin{proof} 
The equality follows from the isomorphism $$\Ext iRMN \cong \Ext {i-1}{R/xR}M{N/xN}$$ for each $i>0$; see \cite[Lemma 2, p.
140]{mat}.
\end{proof}

\begin{proof}
[\rm {\it Proof of Theorem \ref{cx-inequality}}]
First, we show that we may pass to the completion of $R$. Remark that the complexities
involved are not changed by completion. We only need to check
that $\pd {\wh R_P}{\wh M_P} <\infty $ for any prime ideal $P$ in
the punctured spectrum of $\wh R$. Let $p = P\cap R$. Then $p$ is a
prime ideal in the punctured spectrum of $R$, thus $\pd {R_p} {M_p}
<\infty$. But we have the commutative diagram:
\[
\xymatrixrowsep{2pc} \xymatrixrowsep{2pc} \xymatrix{ &R
\ar@{->}[d]\ar@{->}[r] &\wh R\ar@{->}[d] \\ &R_p \ar@{->}[r] & {\wh
R_P}}
\]
Note that the map $R_p \to \wh R_P$ is flat. So $\pd {\wh R_P}{\wh
M_P} <\infty $ as desired.

Since we may assume $R$ is complete,  by the discussion after
\cite[Theorem A]{ab89} there exists a short exact sequence of
$R$-modules $0\to N\to X\to N'\to 0,$ where $X$ is of finite
injective dimension and $N'$ is a maximal Cohen-Macaulay module. By
applying the functor $\Hom RM-$ to this sequence we get from the
long exact sequence of Exts the isomorphism $\Ext iRMN\cong\Ext
{i+1}RM{N'}\quad\text{for all}\ i>\dim R$. Thus without loss of
generalization, we may assume $N$ is a maximal Cohen-Macaulay
module.

Second, by replacing $M$ with a high syzygy  we may assume that $M$
is also a maximal Cohen-Macaulay, and hence $M_p$ is free for all
primes $p\not=\fm$.

Finally, we apply Lemmas \ref{x-reduction} and \ref{plex_deform} to
reduce to the artinian case  and Lemma \ref{cx-px-inequalities} to
finish the proof.
\end{proof}

%%%%%%%%%%%%%%%%%%%%%%%%%%%%%%%%%%%%%%%% SECTION %%%%%%%%%%%%%%%%%%%%%%%%%%%%%%%%%%%%%%%%%%

\section{Equalities of complexities of pairs of modules}
\label{AAR}

In this section $(R, \fm, k)$ is a local ring.

\begin{definition} We define the {\it asymptotic Auslander-Reiten} property of a local ring $R$ to be the following:

\noindent (AAR)
\hspace{0.5cm} For any finitely generated $R$-module
$M$ the implication
\begin{equation*}
\cxx RMR=0=\cxx RMM \implies \cx RM=0
\end{equation*}
\hspace{1.7cm} holds.
\end{definition}

Recall that a local ring $R$ satisfies the Auslandrer-Reiten
condition if it has the following property:

\vspace{0.3cm}

\noindent (AR)
\hspace{0.5cm} For any  finitely generated $R$-module
$M$ the implication
\begin{equation*}
\Ext {i}RMR=0=\Ext iRMM\  \text{for all}\ i>0 \implies M\ \text{is
free}
\end{equation*}
\hspace{1.4cm} holds.

\begin{remark}
\label{AR and AAR}
It is easy to see that if $R$ is a ring satisfying the (AAR)
condition, then it satisfies the (AR) condition. Indeed, let $M$ be
a finitely generated $R$-module with $\Ext {i}RMR=0=\Ext iRMM\
\text{for all}\ i>0$. Since $R$ satisfies (AAR) condition we obtain
that $\pd RM<\infty$. In this case, we know by \cite[(2.6)]{i69}
that $\Ext{\pd RM}RMR\not=0$, thus $M$ is free.

However, we do not know if the reverse implication holds.
\end{remark}

\begin{remark}
\label{AAR-rings}
The Auslander-Reiten Conjecture asserts that every
local ring satisfies (AR). This conjecture has been studied
extensively in the recent papers \cite{abs05,ch,hsv04,lh04}. In
fact, some results in those papers implicitly provide classes of
rings satisfying (AAR). One  such class consists of rings with
radical cube zero; see \cite[(4.1)]{hsv04}. Recently, Christensen
and Holm  show that (AAR) is implied by the Auslander's condition on
the vanishing of cohomology which they denote by (AC); see
\cite[(2.3)]{ch}.
\end{remark}

In this section we investigate rings satisfying stronger properties
than (AAR):
\begin{equation*} \tag{P1} \label{property1} \cxx RMR=\cx
RM\quad\text{for all finitely generated $R$-modules}\  M.
\end{equation*}
\begin{equation*}
\tag{P2} \label{property2} \cxx RMM=\cx RM\quad\text{for all
finitely generated $R$-modules}\  M.
\end{equation*}

\begin{remark}
\label{remark Gorenstein}
A ring with property \eqref{property1} cannot be Gorenstein unless
it is a regular ring. Indeed, assume that $R$ satisfies
\eqref{property1} and is Gorenstein. Therefore, since $\px RR =0$ it
follows that $\pd Rk<\infty,$ thus $R$ is  regular.

On the other hand, if $R$ is a complete Cohen-Macaulay local ring with property
\eqref{property2}, then $R$ is Gorenstein. 
By assumptions $R$ has a canonical module $D$. Since
$D$ has finite injective resolution, we get that $\cxx RDD=0$, thus
$\cx RD=0$. In particular, the module $D$ has finite projective
dimension and finite injective dimension, hence $R$ is Gorenstein by
\cite[(4.4)]{f77}.

However,  there exists an artinian Gorenstein local ring  $R$ not
satisfying \eqref{property2}. If $R$ is artinian, then  we have $\cxx RMM = \cxx
R{M^\vee}{M^\vee}$; see \ref{cx-px}. Therefore, if $R$ satisfies \eqref{property2},
then $\cx RM =\px RM$ for all finitely generated $R$-modules $M$.
However, Jorgensen and \c{S}ega construct in \cite[(1.2)]{js05} a
Gorenstein ring with $\fm^4=0\not=\fm^3$ and a finitely generated
$R$-module $M$ with $1=\cx{R}M<\px{R}{M}=\infty.$ Thus, $R$  does
not satisfy \eqref{property2}.
\end{remark}

In the next two results we identify classes of rings  satisfying the property \eqref{property1}.

\begin{proposition}
\label{2r>l}
Let $R$ be an artinian local ring such that $2\ell_R(\Soc(R))> \ell_R(R)$. Then, for all finitely generated $R$-module $M$
$$ \cxx RMR = \cx RM\in\{0,\infty\}.$$
In particular, if $R$ is not a field, then $\cx RE=\cx Rk=\infty.$
Here $E$ denotes the injective envelope of $k$.
\end{proposition}

\begin{proof}
Set $r=\ell_R(\Soc(R))$ and $l=\ell_R(R)$. It is proved in
\cite[(4.2.7)]{a98} that for every finitely generated  $R$-module
$M$, we have
$$\beta_{i+1}(M) \geq \frac{r}{l-r} \beta_i(M)\quad\text{for all}\quad i>0.$$
In particular, if  $2r > l$ we  obtain
$$\{\cx RM\mid M\ \text{is a finitely generated $R$-module}\} \subseteq \{0,\infty\}.$$
Corollary \ref{cx MN}(1) applied to the pair of modules $(M,R)$, and
the ideal $I=\Soc(R)$ gives  the equality $\cxx RMR=\cx RM$.

If $R$ is a complete intersection artinian ring with  $2r>l$ then
$R$ is a field. If $R$ is not a complete intersection, then $\cx
Rk=\infty$ by \ref{k ci}. By \ref{cx-px} , $\cx RE=\px RR$.
Corollary \ref{cx MN}(2) applied to the pair of modules $(k,R)$, and
$I=\Soc(R)$ gives the equality $\cx RE=\px Rk=\cx Rk=\infty$.
\end{proof}

\begin{remark}
\label{LJ}
One can find examples of local rings $R$ satisfying the hypotheses
of Proposition \ref{2r>l} by taking $R= S/\fn^h$ with $(S,\fn)$ a
regular local ring of Krull dimension at least $3h-4$ and $h\geq 2$.
From the point of view of  Jorgensen and Leuschke \cite{jl07}, the
artinian local rings $R$ with $2\ell_R(\Soc(R))> \ell_R(R)$ are
furthest from being Gorenstein. Indeed, in  \cite[(3.4)]{jl07}, they
define $\mathfrak{g}(R) = \curv RE/ \curv Rk$ where $\curv RM$
denotes the curvature of $M$, a measure of the exponential growth of
the Betti numbers of $M$; see \cite[Ch 4]{a98}. They show that there
are inequalities $0\leq \mathfrak{g}(R) \leq 1 $ and that $R$ is
Gorenstein if and only if $\mathfrak{g}(R)=0$. Thus, the invariant
$\mathfrak{g}(R)$ measures how far $R$ is from being Gorenstein. One
can show, as in the proofs of Propositions \ref{2r>l} and
\ref{sequences}, that $\curv RE=\curv Rk$ when $2\ell_R(\Soc(R))>
\ell_R(R)$. Thus, these rings satisfy $\mathfrak{g}(R) = 1$, and
thus are furthest from being Gorenstein with respect to this
invariant.
\end{remark}

\begin{theorem}
\label{cx m3}
Let $(R,\fm,k)$ be a non-Gorenstein local ring and let $M$ be a finitely generated $R$-module.
\begin{enumerate}[\rm\quad(1)]
\item
If $\fm^2=0$, then
$\cxx RMR = \cx RM \in\{0,\infty\}.$
\item
If $\fm^3=0\not=\fm^2$ and $2\ell_R(\Soc(R))> \ell_R(R)-2$, then
$$\cxx RMR = \cx RM \in\{0,1,\infty\}.$$
\end{enumerate}
\end{theorem}

\begin{proof}
(1): If $M$ is a free module then the statement is clear.
If $M$ is not free, we may assume that $M$ is a finite $k$-vector space
by replacing $M$ by its first syzygy.
Hence, we  obtain the equalities $\cx RM = \cx Rk =\infty$
and $\cxx RMR=\px RR$;  for the second equality we use \ref{k ci}. By hypothesis $R$ is not Gorenstein,
therefore the injective envelope of the residue field $k$ has
infinite complexity by \ref{E Gorenstein}. Thus, by \ref{cx-px} we
have  $\px RR=\infty$, and the desired conclusion follows.

(2): For the rest of the proof, set $b_i=\beta_i^R(M)$ for $i\geq 0$, $r=\ell_R(\Soc(R))$ and $l=\ell_R(R)$. By \cite[Theorem B and (3.9)]{l85} we have
$$\{\cx RM\mid M\ \text{is a finitely generated $R$-module}\}\subseteq \{0,1,\infty\}.$$

 Lemma \ref{cx-px-inequalities} gives the inequality $\cxx RMR\leq \cx RM$. Thus, we may consider the following three cases on complexity of $M$.

If $\cx RM=0$, then $M$ is free, hence $\cxx RMR=0$ by definition.

If $\cx RM=1$, then $\cxx RMR=1$. Otherwise, if  $0=\cxx RMR<\cx RM$ it follows by \cite[(Theorem A)]{cv06}  that  $2r=l-2$, contradicting with our hypothesis.

The last case  is $\cx RM=\infty$. By Proposition \ref{2r>l}, we may assume that
$2r\in\{l,l-1\}$; thus we analyze the two possibilities.

Assume $2r=l$. Lemma \ref{N-inequality} applied to $N=R$ and
$I=\Soc(R)$ gives
\begin{equation*}
\ell_R (\Ext iRMR) \geq  r\cdot (b_i-b_{i-1}).
\end{equation*}
Since the sequence $\{b_i\}_{i\geq 1}$ has infinite complexity, so
does the sequence $\{b_{i+1}-b_i\}_{i\geq 1}$   by Proposition
\ref{sequences}(2).  Thus, $\cxx RMR=\infty$.

Now, assume  $2r=l-1$. Lemma \ref{N-inequality} applied to $N=R$ and
$I=\Soc(R)$ gives
\begin{equation*}
\tag{$*$}
\ell_R(\Ext iRMR) \geq  r\cdot b_i - (r+1)\cdot b_{i-1}.
\end{equation*}
Moreover, by replacing $M$ with its first syzygy, we may assume that $\fm^2M=0$.

Assume that $k$ is not a direct summand of any syzygy of $M$. Set $a=\ell_R(\fm^2)$ and $e=\edim R$.  By \cite[(3.2)]{l85} we have  $r=a$, and by our hypothesis $r=e$.
By \cite[(3.3)]{l85} we know that  the sequence $\{b_i\}_{i\geq 0}$ satisfies
$$b_{i+1}=eb_{i}-ab_{i-1}=r(b_{i}-b_{i-1}),\quad\text{for all}\quad i\geq 1.$$
Therefore, the inequality $(*)$ becomes
\begin{align*}
\ell_R(\Ext iRMN)&\geq  r\cdot b_i - (r+1)\cdot b_{i-1}\\
                &=r\cdot(b_i-b_{i-1})-b_{i-1}\\
                &=b_{i+1}-b_{i-1}\\
                &=(b_{i+1}-b_i)+(b_i-b_{i-1}).
\end{align*}
Thus, we get as above that $\cxx RMR=\infty$.

Finally, assume that for some $j\geq 0$ the $j$-th syzygy of $M$,
denoted $M_j$, satisfies  $M_j\cong k\oplus M_j'$ . Then there are
isomorphisms
$$\Ext iRMR\cong \Ext {i-j}R{M_j}R\cong \Ext {i-j}RkR\oplus\Ext{i-j}R{M_j'}R\quad\text{for all}\quad i>j.$$
Since $R$ is non-Gorenstein with $\fm^3=0$, we have $\cxx RkR=\cx
RE=\infty$; the first equality is by \ref{cx-px} and the second by
\ref{E Gorenstein}. Therefore, $\cxx RMR=\infty$, so we have the
desired conclusion.
\end{proof}

\begin{remark}
\label{JS-example}
The inequality $2\ell_R(\Soc(R))>\ell_R(R)-2$ of
Theorem \ref{cx m3}(2) is sharp. Jorgensen and \c{S}ega construct in
\cite[(3.1)]{js05} a non-Gorenstein ring $R$ with
$\fm^3=0\not=\fm^2$ and $2\ell_R(\Soc(R))=\ell_R(R)-2$ and a
finitely generated $R$-module $M$ with $$0=\cxx RMR<\cx RM=1.$$
\end{remark}

A class of rings satisfying the property \eqref{property2} is identified below.

\begin{proposition}
\label{cx m3-Gorenstein}
Let $(R,\fm,k)$ be a local Gorenstein ring
with  $\fm^3=0$. If $M$ is a finitely generated $R$-module,
then $$\cxx RMM = \cx RM.$$
\end{proposition}

\begin{proof}
If $R$ is a complete intersection, then apply \ref{properties-ci}.
Assume that $R$ is not a complete intersection. Let $N$ be the first
syzygy of $M$, then $\cxx RNN=\cxx RMN$ and $\cx RM=\cx RN$.
Since $\Ext{i}RMR =0$ for all $i>0$, we have $\cxx RMN =\cxx RMM$.
Therefore, it is enough to show that $\cxx RNN = \cx RN$.

If $\fm^2=0$, then $N$ is a $k$-vector space.
Thus, the desired equality follows by the
definition of complexity.

If $\fm^3=0\not=\fm^2$, then $\fm^2N=0$. By Corollary \ref{m2N} we have the inclusion $\cxx RNN\in\{\cx RN, \px RN\}$. But by the discussion in \cite[Sec.2]{js05} we have $\cx RN = \px RN$. Thus, the desired equality holds.
\end{proof}

Combining Theorem \ref{cx m3}(1) and Proposition \ref{cx
m3-Gorenstein} we obtain that over a local ring $(R, \fm)$  with
$\fm^2=0$, every finitely generated $R$-module $M$ satisfies the
equalities  $\cxx RMR=\cx RM=\cxx RMM.$ As we have seen in Remark
\ref{JS-example} this is not true for all the rings with $\fm^3=0$.
For such rings, the (AAR) condition is implied by the following
result, which gives  more information on complexities:

\begin{proposition}
Let $(R, \fm, k)$ be a ring with  $\fm^3=0\not=\fm^2$ and let $M$ be a finitely generated $R$-module such that $\cxx{R}MM=0$, then
$$\cxx{R}MR=\cx{R}M\in\{0,1,\infty\}.$$
\end{proposition}

\begin{proof}
By \cite[Theorem B and (3.9)]{l85}, we have
$\cx RM\in\{0,1,\infty\}.$

If $\cx{R}M=0$, then $\cxx RMR=0$ by definition.

If $\cx{R}M=1$, then $\cxx{R}MR\in\{0,1\}$ as we have $\cxx RMR\leq\cx RM$; see Lemma \ref{cx-px-inequalities}. If $\cxx RMR=0$, then $M$ is free by \cite[(4.1.1)]{hsv04}, contradiction. Thus, in this case we have $\cxx RMR=1$, as desired.

Finally, we consider the case $\cx RM=\infty$. Let $N$ be the first
syzygy of the module $M$; it satisfies $\fm^2N=0$.

If $\ell_R(\fm N) > \nu_R(N)$, then by part (1) of Proposition
\ref{m2-cx} we have the first equality in $\cxx RNN=\cx RN=\cx
RM=\infty$. On the other hand, the long exact sequence obtained  by
applying the functor $\Hom RM-$ to the short exact sequence $0\to
N\to R^{\nu_R(M)}\to M\to 0$ and the assumption $\cxx{R}MM=0$
implies that
$$\Ext iRMN \cong \Ext iRM{R}^{\nu_R(M)} \quad \text{for all} \quad i\gg 0.  $$
It follows that $\cxx{R}MN=\cxx{R}MR$ and this is equal to 
$\cxx{R}NN$; recall that $N$ is a syzygy of $M$. Therefore,  $\cxx{R}MR=\infty$ as desired.

If $\ell_R(\fm N) \leq \nu_R(N)$, then by parts (2) and (3) of
Proposition \ref{m2-cx} we have $\cxx{R}MR=\cxx{R}NR\in\{\px{R}R,
\px{R}R-1\}$. If $R$ is Gorenstein, then by \cite[(4.1.2)]{hsv04} and by hypothesis $M$ is free, contradicting  our assumption. If $R$
is not Gorenstein, then $\px{R}R=\cx RE=\infty$; the first equality
is by \ref{cx-px} and the second  by \ref{E Gorenstein}. Therefore,
$\cxx{R}MR=\infty$.
\end{proof}

\begin{remark}
If $R$ is a complete intersection local ring and $M$  is a finitely
generated $R$-module, then the condition $\cxx RMM=0$ implies by
\ref{properties-ci}  that $\cxx RMR=\cx RM=0$.
\end{remark}

Finally, for completeness, we note that the main arguments in
\cite{ab00} give a slightly more general result than what is stated
in \ref{properties-ci}. Recall that the complete intersection
dimension of a module $M$ is defined as:
$$\ci RM = \inf \{\pd{Q}{M\otimes_RR'}-\pd QR' | \  R \to R' \leftarrow Q\  \text{is a quasi-deformation} \}$$
Here a quasi-deformation $R \to R' \leftarrow Q$ is a diagram of
local homomorphisms such that $R \to R'$ is flat and $R' \leftarrow
Q$ is surjective with kernel generated by a regular $Q$-sequence
$\bsf=f_1,\cdots,f_c$; see \cite{agp97}.

\begin{proposition}
\label{ci}
Let $R$ be a local ring and let $M$ be a finitely generated $R$-module.
If $\ci{R}{M} < \infty$  then $\cxx RMM=\cx RM$.
\end{proposition}

\begin{proof}
By definition \cite[Section 4]{ab00}, there exists a quasi-deformation
as above such that  $\pd{Q}{M'}<\infty$, where $M'=M\otimes_R R'$.
Replacing $R,M$ by $R', M' $ we may assume $R=R'$. We may also
assume that $k$ is algebraically closed by replacing $R$ by its
residual algebraic closure; see \cite[(4.1.1)]{ab00} and
\cite[(1.14)]{agp97}.

Now, by  \cite[(2.4)]{ab00}, we have for any $R$-module
$N$ the equality
$$\cxx RMN = \dim V^*(Q,\bsf,M,N).$$
Here $  V^*(Q,\bsf,M,N)$ denotes the support variety of the pair
$(M,N)$; see \cite[(2.1)]{ab00}. By  \cite[(2.5)]{ab00}, one has an
equivalent definition: $$ V^*(Q,\bsf,M,N) = \{ a\in {k}^r \Big|
\Ext{n}{Q_a}{M}{N} \neq 0 \ \text{for infinitely many}\ n \} \cup \{
0\},$$ where
  $a=( a_1,\cdots, a_r), f_a = \sum a_if_i$ and $Q_a = {Q}/(f_a)$.
To prove $\cxx RMM  =\cx RM$ it suffices to show $V^*(Q,\bsf,M,M) =
V^*(Q,\bsf,M,k)$. By the above definition we have to show that for each
$a\in {k}^r$,  $\Ext{n}{Q_a}{M}{M} =0 \ \text{for} \ n\gg 0
\Leftrightarrow \pd {Q_a}M <\infty$. But this follows directly from
\cite[(4.2)]{ab00} which asserts that a finite module $M$ of finite
CI-dimension over a Noetherian ring $R$ has finite projective
dimension if and only if $\Ext{2i}RMM=0$ for some $i>0$; note that
$\ci{Q_a}{M}<\infty$ by definition.
\end{proof}

\section*{Acknowledgments}
The authors would like to thank Melvin Hochster for his help in
proving Lemma \ref{same h}. We are also grateful to Lars Winther
Christensen for many helpful comments on a preliminary version of
the manuscript. We would like to thank the anonymous referees
for helpful and very detailed comments which improved the paper.

\end{document}